\documentclass[11pt]{amsart}
\usepackage{amscd}
\usepackage{amsmath}
\usepackage{graphicx}
\usepackage{amsfonts}
\usepackage{amssymb}
\begin{document}
\title
{Samuel multiplicities and  Browder Spectrum of  Operator Matrices}

\author{Shifang Zhang}
\address[Shifang Zhang]{
Department of Mathematics, Zhejiang University, Hangzhou 310027, P.R.
China\\School of Mathematics and Computer Science, Fujian Normal University, Fuzhou 350007, China} \email{shifangzhangfj@163.com}

\author{Junde Wu}
\address[Junde Wu]{Department of Mathematics, Zhejiang University, Hangzhou 310027, P R. China}
\email{wjd@zju.edu.cn}

\thanks{{\it This work is
supported by the NSF of China (Grant Nos. 10771034 and 10771191)}}

\begin{abstract} In this paper, we first point out that the necessity
of Theorem 4 in [8] does not hold under the given condition
and present a revised version with a little modification. Then we show
that the definitions of some classes of semi-Fredholm operators,
which use the language of algebra and first introduced by X. Fang in
[8], are equivalent  to  that of some well-known operator
classes. For example, the concept of shift-like semi-Fredholm
operator on Hilbert space coincide with that of upper semi-Browder
operator. For applications of Samuel
multiplicities we characterize the sets of $\bigcap_{C\in
B(K,\,H)}\sigma_{ab}(M_{C}),\bigcap_{C\in
 B(K,\,H)}\sigma_{sb}(M_{C})$ and $\bigcap_{C\in
 B(K,\,H)}\sigma_{b}(M_{C}),$ respectively, where
$M_{C}=\left(\begin{array}{cc}A&C\\0&B\\\end{array} \right)$ denotes
a 2-by-2 upper triangular operator matrix acting on the Hilbert
space $H\oplus K$.

\end{abstract}

\maketitle
\section{Introduction}
Throughout this paper, let $H$ and $K$ be separable infinite
dimensional complex Hilbert spaces and $B(H, K)$ the set of all
bounded linear operators from $H$ into $K$, when $H=K$, we write
$B(H, H)$ as $B(H)$. For $A\in B(H)$, $B\in B(K)$ and   $C\in B
(K, H)$, we have $M_{C}=\left(\begin{array}{cc}A&C\\0&B\\\end{array} \right)\in B(H\oplus K)$.
 For $T\in B(H, K)$, let $R(T)$ and $N(T)$ denote the range and kernel of $T$, respectively, and denote $\alpha
(T)=\dim N(T)$, $\beta (T)=\dim K / R(T)$. If $T\in B(H)$, the
ascent $asc(T)$ of $T$ is defined to be the smallest nonnegative
integer $k$ which satisfies that $N(T^{k})=N(T^{k+1})$. If such $k$
does not exist, then the ascent of $T$ is defined as infinity.
Similarly, the descent $des(T)$ of $T$ is defined as the smallest
nonnegative integer $k$ for which $R(T^{k})=R(T^{k+1})$ holds. If
such $k$ does not exist, then $des(T)$ is defined as infinity, too.
If the ascent and the descent of $T$ are finite, then they are equal
(see [3]). For $T\in B(H)$, if $R(T)$ is closed and
$\alpha (T)<\infty$, then $T$ is said to be a upper semi-Fredholm
operator, if $\beta (T)<\infty$, which implies that $R(T)$ is
closed,  then $T$ is said to be a lower semi-Fredholm operator. If
$T\in B(H)$ is either upper or lower semi-Fredholm operator, then
$T$  is said to be a semi-Fredholm operator. If both $\alpha
(T)<\infty$ and $\beta (T)<\infty$, then $T$ is said to be a
Fredholm operator. For a semi-Fredholm operator $T$, its index ind
$(T)$ is defined by ind $(T)=\alpha(T )-\beta(T).$

In this paper, the sets of invertible operators, left invertible
operators and right invertible operators on $H$ are denoted by
$G(H), G_l(H)$ and $G_r(H)$, respectively, the sets of all Fredholm
operators, upper semi-Fredholm operators and lower semi-Fredholm
operators on $H$ are denoted by $\Phi(H)$, $\Phi_{+}(H)$ and
$\Phi_{-}(H)$, respectively, the sets of all Browder operators,
upper semi-Browder operators and lower semi-Browder operators on $H$
are defined, respectively, by
\begin{eqnarray*}
& & \Phi_{b}(H):=\{T \in  \Phi(H): asc(T)= des(T)<\infty \},\\
& &  \Phi_{ab}(H):=\{T \in  \Phi_{+}(H): asc(T)<\infty \},\\
& &\Phi_{sb}(H):=\{T \in  \Phi_{-}(H): des(T)<\infty \}.
\end{eqnarray*}
Moreover, for $T\in B(H)$, we introduce its corresponding spectra as
following [19]:

\vspace{3mm}

\par \par the spectrum: $\sigma_{}(T)=\{\lambda \in
{\mathbb{C}}:T-\lambda I\not \in G(H) \}$,
\par the left spectrum: $\sigma_{l}(T)=\{\lambda \in {\mathbb{C}}:T-\lambda I\not \in G_l(H)\}$,
\par the right spectrum: $\sigma_{r}(T)=\{\lambda \in {\mathbb{C}}:T-\lambda I\not \in G_r(H)\}$,
\par the essential spectrum:   $\sigma_{e}(T)=\{\lambda \in {\mathbb{C}}:T-\lambda I\not \in \Phi(H)\}$,
\par the upper semi-Fredholm spectrum: $\sigma_{SF+}(T)=\{\lambda \in {\mathbb{C}}: T-\lambda I \not \in \Phi_+(X)\},$
\par the lower semi-Fredholm spectrum:  $\sigma_{SF-}(T)=\{\lambda \in {\mathbb{C}}: T-\lambda I \not \in \Phi_-(X)\},$
\par the Browder spectrum: $\sigma_{b}(T)=\{\lambda \in {\mathbb{C}}:T-\lambda I \not \in \Phi_{b}(H)\},$
\par the upper semi-Browder spectrum: $\sigma_{ab}(T)=\{\lambda \in {\mathbb{C}} : T-\lambda I \not \in \Phi_{ab}(X)\},$
\par the lower semi-Browder spectrum: $\sigma_{sb}(T)=\{\lambda \in {\mathbb{C}} : T-\lambda I \not \in \Phi_{sb}(X)\}.$

\vskip 0.2in

For a semi-Fredholm operator $T\in B(H)$, its shift Samuel
multiplicity $s\_{mul}(T)$ and backward shift Samuel multiplicity
$b.s\_{mul}(T)$ are defined ([5-8]),
respectively, by
$$s\_{mul}(T)=\lim\limits_{k \rightarrow \infty}\frac{\beta(T^k )}{k},$$
$$b.s\_{mul}(T)= \lim\limits_{k\rightarrow \infty}\frac{\alpha(T^k
)}{k}.$$ Moreover, it has been proved that $s\_{mul}(T),
b.s.\_{mul}(T)\in \{0,1,2,\ldots,\infty\} \makebox{\,\,\,and\,\,\,
ind}(T)= b.s.\_{mul}(T)-s\_{mul}(T)$. These two invariants refine
the Fredholm index and can be regarded as the stabilized dimension
of the kernel and cokernel [8].

\vskip 0.2in

\noindent {\bf Definition 1.1} ([8]). A semi-Fredholm
operator $T\in B(H)$ is called a pure shift semi- Fredholm operator
if $T$ has the form $T=U^nP$, where $n\in \bf N$ or $n=\infty$,  $U$
is the unilateral shift, and $P$ is a positive invertible operator.
Analogously, $T$ is called a pure backward shift semi-Fredholm
operator if its adjoint $T^*$ is a pure shift semi-Fredholm
operator. Here $U^\infty$ denotes the direct sum of countably (infinite) many copies of  $U$.

\vskip 0.2in \noindent {\bf Definition 1.2} ([8]) A
semi-Fredholm operator $T\in B(H)$ is called a shift-like
semi-Fredholm operator if $b.s.\_{ mul}(T)=0$; $T$ is called a shift
semi-Fredholm operator if $N (T)=0$. Analogous concepts for backward
shifts can also be defined. $T$ is called a stationary semi-Fredholm
operator if $b.s.\_{ mul}(T)=0$ and $s\_{ mul}(T)=0$.

\vskip 0.2in

It follows from Definition 1.1 that $T$ is a shift semi-Fredholm
operator iff $T$ is a left invertible operator, and that $T$ is a
backward shift semi-Fredholm operator iff $T$ is a right invertible
operator.

\vskip 0.2in

In ([8], Theorem 4 and  Corollary 18), Fang gave the following
$4\times 4$ upper-triangular representation theorem: An operator
$T\in B(H)$ is semi-Fredholm iff $T$ can be decomposed into the
following  form with respect to some orthogonal decomposition
$H=H_1\oplus H_2\oplus H_3\oplus H_4,$
$$T=\left(\begin{array}{cccc}T_1&*&*&*\\0&T_2&*&*\\0&0&T_3&*\\0&0&0&T_4\end{array}\right),$$ where  $\dim H_4<\infty$, $T_1$ is a pure backward shift
semi-Fredholm operator, $T_2$ is invertible, $T_3$ is a pure shift
semi-Fredholm operator, $T_4$ is a finite nilpotent operator.
Moreover, ind$(T_1)= b.s.\_{mul}(T)$ and ind$(T_3)= -s\_{mul}(T)$.

\vskip 0.2in

The following example shows that the representation theorem is not
accurate.

\vskip 0.2in

\noindent {\bf Example 1.3.} Let $H$ be the direct sum of countably
many copies of $\ell^2:=\ell^2(\bf N)$, that is, the elements of $H$
are the sequences $\{x_j\}_{j=1}^\infty$ with $x_j \in \ell^2$ and
$\sum_{j=1}^\infty \| x_j\|^2 < \infty$. Let $V$ be the unilateral
shift on $\ell^2$, i.e., $$V : \ell^2 \to \ell^2, \quad \{z_1,z_2,
\ldots\} \mapsto \{0,z_1,z_2,\ldots\},$$ and the operators $T_1$ and
$T_3$ be defined by
$$ T_1: H \to H, \quad \{x_1, x_2, \ldots\} \mapsto
\{V^*x_1, V^*x_2, \ldots\}$$
and
$$ T_3: H \to H, \quad \{x_1, x_2,
\ldots\} \mapsto \{Vx_1, Vx_2, \ldots\}.$$ Now, we consider the
operator
$$T=\left(\begin{array}{ll} T_1&0\\0&T_3\end{array}\right):
H \oplus H \to H \oplus H.$$ Note that $T_1$ is a pure backward
shift semi-Fredholm operator, $T_3$ is a pure shift semi-Fredholm
operator, so $T$ satisfies the conditions of Fang's $4\times 4$
triangular representation theorem, but, since
$\alpha(T_1)=\alpha(T)=\beta(T)=\dim(H/ R(T_3))=\infty$, so $T$ is
not a semi-Fredholm operator.

\vskip 0.2in

Now, we can prove the following improved $4\times 4$
upper-triangular representation theorem:

\vskip 0.2in

\noindent {\bf Theorem 1.4.} An operator $T\in B(H)$ is
semi-Fredholm iff $T$ can  be decomposed into the following form
with respect to some orthogonal decomposition $H=H_1\oplus H_2\oplus
H_3\oplus H_4$,
$$T=\left(\begin{array}{cccc}T_1&*&*&*\\0&T_2&*&*\\0&0&T_3&*\\0&0&0&T_4\end{array}\right),$$ where  $\dim H_4<\infty$, $T_1$ is a pure backward shift
semi-Fredholm operator, $T_2$ is invertible, $T_3$ is a pure shift
semi-Fredholm operator and $\min \{\makebox{ind}(T_1),
-\makebox{ind}(T_3)\}<\infty,$ $T_4$ is a finite nilpotent operator.
Moreover,

\vskip 0.1in

\noindent (1) $ \ $ ind$(T_1)=b.s.\_{mul}(T)$, ind$(T_3)=
-s\_{mul}(T)$;

\noindent (2) $ \ $ ind$(T)=+\infty$ iff ind$(T_1)=+\infty$;

\noindent (3) $ \ $ ind$(T)=-\infty$ iff ind$(T_3)=-\infty$;

\noindent (4) $ \ $ ind$(T)$ is finite iff both of ind$(T_1)$ and
ind$(T_3)$ are finite.

\vskip 0.2in

Theorem 1.4 can be described as  $3\times 3$ triangular representation
form which may be more convenient for the study of operator theory,
that is,

\vskip 0.2in

\vskip 0.2in

\noindent {\bf Theorem 1.5}. An operator $T\in B(H)$ is
semi-Fredholm  if and only if $T$ can  be decomposed into the
following  form with respect  to  some orthogonal decomposition
$H=H_1\oplus H_2\oplus H_3$

$$T=\left(\begin{array}{ccc}T_1&T_{12}&T_{13}\\0&T_2&T_{23}\\0&0&T_3\end{array}\right):
H_1\oplus H_2\oplus H_3\rightarrow H_1\oplus H_2\oplus H_3,$$

\noindent where  $\dim H_3<\infty$, $T_1$ is  a right invertible
operator, $T_3$ is a finite, nilpotent operator, $T_2$ is a left
invertible operator, and $\min \{\makebox{ind}(T_1),
-\makebox{ind}(T_2) \}<\infty.$ Moreover,
ind$(T_1)=\alpha(T_1)=b.s.\_mul(T)$, ind$(T_2)=
-\beta(T_2)=-s\_mul(T)$ and ind$(T)= \alpha(T_1) -\beta(T_2)$.

\vskip 0.1in

The next lemma is useful for the proofs of our results below,
especially in Section 2.

\vskip 0.1in

\medskip\noindent {\bf Lemma 1.6} [19].
Let $A\in B(H)$, $B\in B(K)$ and $C\in B(K,H)$.

\noindent (1) $ \ $ If $A\in \Phi_{b}(H)$, then $B\in \Phi_{ab}(K)$
iff
          $M_C\in \Phi_{ab}(H\oplus K)$ for some $C\in B(K, H)$.

\noindent (2) $ \ $ If  $M_C\in \Phi_{ab}(H\oplus K)$ for some $C\in
B(K, H)$,
          then $A\in\Phi_{ab}(H)$.

\noindent (3) $ \ $ If $A\in \Phi_{ab}(H)$ and $B\in \Phi_{ab}(K)$,
then $M_C\in
          \Phi_{ab}(H\oplus K)$ for any $C\in B(K, H)$.

\noindent (4) $ \ $ If $B\in \Phi_{b}(K)$, then $A\in \Phi_{ab}(H)$
iff $M_C\in \Phi_{ab}(H\oplus K)$ for some $C\in B(K, H)$;

$A\in \Phi_{sb}(H)$) iff $M_C\in \Phi_{sb}(H\oplus K)$ for some
$C\in B(K, H)$.

\noindent (5) $ \ $ If  $M_C\in \Phi_{b}(H\oplus K)$ for some $C\in
B(K, H)$,
          then  $A\in \Phi_{ab}(H)$ and $B\in \Phi_{sb}(K)$.

\noindent (6) $ \ $ If two of $A$, $B$ and $M_C$ are Browder, then
so is the third.

\vskip 0.2in \noindent {\bf Proposition 1.7.} Let $T\in B(H)$. Then
$T$ is upper semi-Browder iff $T$ can be decomposed into the
following form with respect to some orthogonal decomposition
$H=H_1\oplus H_2$,
$$T=\left(\begin{array}{cc}T_1&T_{12}\\0&T_2\end{array}\right),$$ where $\dim(H_1)<\infty$, $T_1$ is nilpotent, $T_2$ is
 left invertible, and $\beta(T_2)=s\_mul(T)$=$-$ind$(T)$.

 \vskip 0.2in

 \noindent {\bf Proof.} Necessity. Suppose that $T$ is upper semi-Browder. Then we can assume
$p=asc(T)<\infty$. Let $H_1= N(T^p)$. Note that $T$ is upper
semi-Fredholm, so $\dim H_1<\infty$. Let $H=H_1\oplus H_1^\bot$, we
have
$$T=\left( \begin{array}{cc}  T_1&T_{12}\\ 0&T_2\\ \end{array}
\right):H_1\oplus H_1^\bot \rightarrow  H_1\oplus H_1^\bot.$$ That
$T_1$ is nilpotent is clear. Moreover, since the fact that $\dim
H_1<\infty$ implies $T_1\in\Phi_{b}(H_1)$, it follows from Lemma 1.6
(1) that $T_2\in\Phi_{ab}(H_1^\bot)$. A direct calculation shows
that $T_2$ is injective, thus, $T_2$ is left invertible. From
Theorem 1.5, it is clear that
$\beta(T_2)=s\_mul(T)=\makebox{ind}(T_2)$.

Sufficiency follows from Lemma 1.6 immediately.

\vskip 0.2in

\noindent {\bf Proposition 1.8.} Let $T\in B(H)$. Then $T$ is lower
semi-Browder iff $T$ can be decomposed into the following form with
respect to some orthogonal decomposition $H=H_1\oplus H_2,$
$$T=\left(\begin{array}{cc}T_1&T_{12}\\0&T_2\end{array}\right),$$ where $\dim(H_2)<\infty$, $T_1$ is right invertible,
$T_2$ is nilpotent, and $\alpha(T_1)=b.s.\_mul(T)$=ind$(T)$.

\vskip 0.2in

\noindent {\bf Proof.} Necessity. If $T$ is lower semi-Browder,
then we can assume $p=des(T)<\infty$. Denote $H_1= R(T^p)$ and
$H_2=H_1^\bot$. Note that $T^p$ is lower semi-Browder, so $\dim
H_2<\infty$. Let $H=H_1\oplus H_2$, we have
$$T=\left( \begin{array}{cc}   T_1&T_{12}\\ 0&T_2\\ \end{array}
\right):H_1\oplus H_2 \rightarrow H_1\oplus H_2.$$ That $T_1$ is
surjective and $T_2^P=0$ is evident. Note that $\dim H_2<\infty$
implies $T_2\in\Phi_{b}(H_2)$, it follows from Lemma 1.6 that
$T_1\in\Phi_{sb}( H_1),$ and so $T_1$ is right invertible. From
Theorem 1.5, we have $\alpha(T_1)=\makebox{ind}(T_1)=b.s.\_mul(T)$.

Sufficiency follows from Lemma 1.6.

\vskip 0.1in

Combining Theorem 1.5, Propositions 1.7 and 1.8, we have
the following theorem immediately.

\vskip 0.2in

\noindent {\bf Theorem 1.9.} Let $T\in B(H)$. Then

\noindent (1) $ \ $ $T$ is a shift-like semi-Fredholm operator iff
$T$ is an upper semi-Browder operator.

\noindent (2) $ \ $ $T$ is a backward shift-like semi-Fredholm
operator iff $T$ is a lower semi-Browder operator.

\noindent (3) $ \ $ $T$ is a stationary semi-Fredholm operator iff
$T$ is a Browder operator.

\section{Applications of Samuel multiplicities}

 In ([8-12]), Fang studied Samuel
 multiplicities and presented some applications. In this section,
 by using Samuel multiplicities, we characterize the sets
$\bigcap_{C\in B(K,\,H)}\sigma_{ab}(M_{C})$, $\bigcap_{C\in
B(K,\,H)}\sigma_{sb}(M_{C})$ and $\bigcap_{C\in
B(K,\,H)}\sigma_{b}(M_{C})$ completely, where
$M_{C}=\left(\begin{array}{cc}A&C\\0&B\\\end{array} \right)$ is a
$2\times 2$ upper triangular operator matrix defined on $H\oplus K$.
For the study advances of $2\times 2$ upper triangular operator
matrix, see ([1-4],
[13-19]).

\vskip 0.2in

First, note that if $T\in B(H)$, then $T$ is bounded below iff $T$
is left invertible, thus, Theorem 1 of [14] can be
rewritten as follows:

 \vskip 0.2in

 \noindent {\bf Lemma 2.1}
[14]. For any given $A\in B(H)$ and $ B\in B(K)$, $M_C$ is
left invertible for some $C\in B(K, H)$ iff $A$ is left invertible
and
$\left\{\begin{array}{cc}a(B)\leq\beta(A)&\makebox{\,if\,}R(B)\makebox{\,is
closed,\,} \\ \beta(A)=\infty&\makebox{\,if\,}R(B)\makebox{\,is not
closed.\,}\\\end{array} \right. $

  \vskip 0.2in

\noindent {\bf Lemma 2.2} [4]. For any given $A\in B(H)$
and $ B\in B(K)$,\begin{equation}\bigcap_{C\in
B(K,\,H)}\sigma(M_{C})
=\sigma_{l}(A)\cup\sigma_{r}(B)\cup\{\lambda\in{\mathbb{C}}:
\alpha{(B-\lambda)}\not=\beta{(A-\lambda)}\}.\end{equation}

\vskip 0.2in

One of the main results in this section is:

\vskip 0.2in

\vskip 0.2in \noindent {\bf Theorem 2.3.}${\label{2}}$ For any given
$A\in B(H)$ and $ B\in B(K)$, $M_C\in \Phi_{ab}(H\oplus K)$ for some
$C\in B(K, H)$ iff $A\in \Phi_{ab}(H)$ and

$$\left\{\begin{array}{ll}s\_{mul}(A)=\infty&\makebox{\,\,\,if\,}B\not
\in \Phi_{+}(K), \\b.s.\_{mul}(B)\leq
s\_{mul}(A)&\makebox{\,\,\,if\,}
 B\in \Phi_{+}(K).\\\end{array} \right. $$

\noindent{\bf Proof.} We first claim  that if $B\not \in
\Phi_{+}(K)$, then \begin{equation}M_C\in \Phi_{ab}(H\oplus K) \,\,
\makebox{ for some}\,\, C\in B(K, H)\Leftrightarrow A\in
\Phi_{ab}(H)\,\, \makebox{and}\,\, s\_{mul}(A)=\infty.\end{equation}

\noindent To do this, suppose $M_C\in \Phi_{ab}(H\oplus K)$.  Then
from Lemma 1.6 we have $A\in \Phi_{ab}(H)$. If $s\_{mul}(A)<\infty,$
then $A\in \Phi_{}(H)$, since
ind$(A)=\alpha(A)-\beta(A)=b.s.\_{mul}(A)- s\_{mul}(A)$. Hence it is
easy to show that $B\in \Phi_{+}(K)$, which is in a contradiction.
Thus, $s\_{mul}(A)=\infty.$

Conversely, suppose that $A\in \Phi_{ab}(H)$ and
$s\_{mul}(A)=\infty,$ which implies $\beta(A)=\infty.$ It follows
from Proposition 1.7 that $A$ can be decomposed into the following
form with respect to some orthogonal decomposition $H=H_1\oplus H_2$

$$A=\left(\begin{array}{cc}A_1&A_{12}\\0&A_2\end{array}\right),$$

\noindent where  $\dim(H_1)<\infty$, $A_1$ is nilpotent, and $A_2$
is a left invertible operator. Noting that $\beta(A)=\infty$, we
have $\beta(A_2)=\infty.$  Hence it follows from Lemma 2.1 that
there exists some  $C_0\in B(K, H_2)$ such that
$\left(\begin{array}{cc}A_2&C_0\\0&B\end{array} \right)$ is left
invertible. Now consider operator

$$M_C=\left(\begin{array}{cc} A&C\\0&B\end{array}\right)=\left(
\begin{array}{ccc} A_1&A_{12}&0\\0&A_2&C_0\\0&0&B\end{array}\right),$$
where $C=\left(\begin{array}{c}0\\C_0\\\end{array} \right)\in
B(K,H).$ By Lemma 1.6, it is easy to check that $M_C
\in\Phi_{ab}(H\oplus K).$

\noindent Next, We claim  that if $B\in \Phi_{+}(K)$, then

\begin{equation}M_C\in \Phi_{ab}(H\oplus K) \,\, \makebox{ for some}\,\, C\in
B(K, H)\Leftrightarrow A\in \Phi_{ab}(H)\,\, \makebox{and}\,\,
b.s.\_{mul}(B)\leq s\_{mul}(A).\end{equation}

\noindent To this end, suppose $M_C\in \Phi_{ab}(H\oplus K)$, which
implies $A\in \Phi_{ab}(H)$. By Proposition 1.8, we have that $A$
can be decomposed into the following form with respect to some
orthogonal decomposition $H=H_1\oplus H_2$

$$A=\left(\begin{array}{cc}A_1&A_{12}\\0&A_2\end{array}\right),$$

\noindent where  $\dim(H_1)<\infty$, $A_1$ is nilpotent, $A_2$ is a
left invertible operator, and $\beta(A_2)=s_\_mul(A)$. Since the
assumption that $B\in \Phi_{+}(K)$, using Theorem 1.5, we know that
$B$ can be decomposed into the following  form with respect  to some
orthogonal decomposition $K=K_1\oplus K_2\oplus K_3$

$$B=\left(\begin{array}{ccc}B_1&*&*\\0&B_2&*\\0&0&B_3\end{array}\right),$$

\noindent where  $\dim K_3<\infty$, $B_1$ is  a right invertible
operator, $B_2$ is a left invertible operator, $B_3$ is a finite,
nilpotent operator, and the parts marked by $*$ can be any
operators.  Moreover, ind$(B_1)= \alpha(B_1)=b.s.\_mul(B)$,
$ind(B_2)= -\beta(B_2)=-s\_mul(B_1)$ and ind$(B)= \alpha(B_1)
-\beta(B_2)$. Therefore, $M_C$ can be rewritten as the following
form

$$M_C=\left(\begin{array}{ccccc} A_1&A_{12}&C_{11}&C_{12}&C_{13}\\
 0&A_2&C_{21}&C_{32}&C_{23}\\0&0&B_1&*&*\\0&0&0&B_2&*\\0&0&0&0&B_3
 \end{array}\right):H_1\oplus H_2\oplus K_1\oplus K_2\oplus K_3
 \rightarrow  H_1\oplus H_2\oplus K_1\oplus K_2\oplus K_3.$$

\noindent Noting that $\dim(H_1)<\infty$ and $\dim(K_3)<\infty$, we
have $A_1\in\Phi_{b}(H_1)$ and $B_3\in\Phi_{b}(K_3)$. Consequently,
Lemma 1.6 leads to

$$\left(\begin{array}{ccc}A_2&C_{21}&C_{32}\\0&B_1&*\\
0&0&B_2\end{array}\right)\in\Phi_{ab}(H_2\oplus K_1\oplus K_2),$$

\noindent which implies $$\left(\begin{array}{cc} A_2&C_{21}
\\0&B_1\end{array}\right)\in \Phi_{ab}(H_2\oplus K_1).$$
Now we shall prove that $$\beta(A_2)\geq\alpha (B_1).$$ If
$\beta(A_2)=\infty$, the above inequality obviously holds. On the
other hand, if $\beta(A_2)<\infty$, then $A_2\in \Phi(H_2)$, and
hence $B_1\in \Phi_+(K_1)$. Thus,

$$0\geq \,\,\makebox{ind}(\left(\begin{array}{cc} A_2&C_{21}
\\0&B_1\end{array}\right))= \,\,\makebox{ind} (A_2)+\,\,\makebox{ind}
(B_1)=-\beta(A_2)+\alpha(B_1),$$ that is,
$$\alpha(B_1)\leq\beta(A_2).$$ Therefore,
$$b.s.\_{mul}(B)\leq s\_{mul}(A).$$

Conversely, suppose $A\in \Phi_{ab}(H)$, $B\in \Phi_{+}(K)$ and
$b.s.\_{mul}(B)\leq s\_{mul}(A)$. Similar to the above arguments, we
have
$$A=\left(\begin{array}{cc}A_1&A_{12}\\0&A_2\end{array}\right):H_1\oplus H_2\mapsto H_1\oplus H_2$$
and
$$B=\left(\begin{array}{ccc}B_1&*&*\\0&B_2&*\\0&0&B_3\end{array}\right)
:K_1\oplus K_2\oplus K_3\mapsto K_1\oplus K_2\oplus K_3,$$

\noindent where $\dim(H_1)<\infty$, $A_1$ is nilpotent, $A_2$ is a
left invertible operator; $\dim K_3<\infty$, $B_1$ is a right
invertible operator, $B_2$ is a left invertible operator, $B_3$ is a
finite, nilpotent operator, and the parts marked by $*$ can be any
operators. Moreover, $\beta(A_2)=s\_mul(A)$ and
$\alpha(B_1)=b.s.\_mul(B)$. Since the assumption that
$b.s.\_{mul}(B)\leq s\_{mul}(A)$,  we have $\alpha(B_1)\leq
\beta(A_2)$. It follows from Lemma 2.1 that there exists a left
invertible operator $\widetilde{C} \in B(K_1, H_2)$ such that

$$\left(\begin{array}{cc}A_2&\widetilde{C}\\0&B_1\end{array}\right)\in
B(H_2\oplus K_1)\,\, \makebox{is left invertible.}$$

\noindent Consider operator $M_C=\left(\begin{array}{cc} A&C\\
 0&B \end{array}\right):H\oplus K \rightarrow H\oplus K$
$$\,\,\,\,\,\,\, \qquad=\left(\begin{array}{ccccc} A_1&A_{12}&0&0&0\\
 0&A_2&\widetilde{C}&0&0\\0&0&B_1&*&*\\0&0&0&B_2&*\\0&0&0&0&B_3
 \end{array}\right):H_1\oplus H_2\oplus K_1\oplus K_2\oplus K_3
 \rightarrow  H_1\oplus H_2\oplus K_1\oplus K_2\oplus K_3,$$

\noindent where $C=\left(\begin{array}{ccc}
0&0&0\\\widetilde{C}&0&0\end{array}\right)\in B(K_1\oplus K_2\oplus
K_3,H_1\oplus H_2).$  Using Lemma 1.6 , it is easy to see that $M_C
\in \Phi_{ab}(H\oplus K).$

\vskip 0.1in By duality, we have

\vskip 0.2in \noindent {\bf Theorem 2.4}${\label{2}}$. For any given
$A\in B(H)$ and $ B\in B(K)$, $M_C\in \Phi_{sb}(H\oplus K)$ for some
$C\in B(K, H)$ iff $B\in \Phi_{sb}(K)$ and
$$\left\{\begin{array}{ll}b.s.\_mul(B)=\infty&\makebox{\,\,\,if\,}
A\not \in \Phi_{-}(H) \\ b.s.\_mul(B)\geq s\_mul(A)&\makebox{\,\,\,
if\,}A\in \Phi_{-}(H)\\\end{array} \right. $$

\vskip 0.1in From Theorems 2.3 and 2.4, we obtain the following two
corollaries, concerning perturbations of the upper semi-Browder
spectrum and  lower semi-Browder spectrum, respectively.

 \vskip 0.2in

\noindent {\bf Corollary 2.5. } For any given $A\in B(H)$ and $ B\in
B(K)$, we have

$$\bigcap_{C\in B(K,\,H)}\sigma_{ab}(M_{C})=\sigma_{ab} (A)\cup
\{\lambda\in {\mathbb{C}}:\lambda\in\sigma_{SF+}(B)\,\makebox{and}~
s.\_mul(A-\lambda)<\infty \} \cup$$$$\{\lambda\in\Phi(A)\cap\Phi_+(B):
 b.s.\_mul(B-\lambda)>s.\_mul(A-\lambda)\}.$$

\vskip 0.2in \noindent {\bf Corollary 2.6. } For any given $A\in
B(H)$ and $ B\in B(K)$, we have

$$\bigcap_{C\in B(K,\,H)}\sigma_{sb}(M_{C})=\sigma_{sb} (B)\cup
\{\lambda\in{\mathbb{C}}:\lambda \in\sigma_{SF-}(A)\,\makebox{and}\,
b.s.\_mul(B-\lambda)<\infty \} \cup$$$$\{\lambda\in\Phi(B)\cap
\Phi_-(A):b.s.\_mul(B-\lambda)<s.\_mul(A-\lambda)\}.$$

 \vskip 0.2in

\noindent {\bf Theorem 2.7.} For any given $A\in B(H)$ and $ B\in
B(K)$, the following statements are equivalent:

 \vskip 0.2in

\noindent (1) $ \ $ $M_C\in \Phi_{b}(H\oplus K)$ for some $C\in B(K,
H)$;

\noindent (2) $ \ $  $A\in \Phi_{ab}(H)$, $B\in \Phi_{sb}(K)$ and
     $b.s.\_mul(B)=s\_mul(A)$;

\noindent (3) $ \ $ $A\in \Phi_{ab}(H)$, $B\in \Phi_{sb}(K)$ and
      $\alpha(A)+ \alpha(B)= \beta(A)+\beta(B)$.

 \vskip 0.2in

\noindent{\bf Proof.} $(1)\Rightarrow(2)$. Suppose that $M_C\in
\Phi_{b}(H\oplus K)$. Then from Lemma 1.6, we have
$A\in\Phi_{ab}(H)$ and $B\in \Phi_{sb}(K)$. Using Propositions 1.7
and 1.8, we have

$$M_C=\left(\begin{array}{cccc} A_1&A_{12}&C_{11}&C_{12}\\
 0&A_2&C_{21}&C_{32}\\0&0&B_1&B_{12}\\0&0&0&B_2
 \end{array}\right):H_1\oplus H_2\oplus K_1\oplus K_2
 \rightarrow  H_1\oplus H_2\oplus K_1\oplus K_2,$$

\noindent where $\dim(H_1)<\infty$, $A_1$ is nilpotent, $A_2$ is a
left invertible operator, $\dim K_2<\infty$, $B_1$ is a right
invertible operator, $B_2$ is a finite, nilpotent operator.
Moreover, $$\beta(A_2)=s.\_mul(A) \,\,\makebox{and}\,\,
\alpha(B_1)=b.s.\_mul(B).$$

\noindent In addition, it follows from Lemma 1.6 that

$$\left(\begin{array}{cc} A_2&C_{21}\\0&B_1\\
\end{array}\right)\in \Phi_{b} (H_2\oplus K_1).$$
Note the well-known fact that if $M_C\in\Phi_{} (H\oplus K)$, then
$A\in\Phi_{} (H)$ if and only if $B \in\Phi(K)$. Thus, if
$\beta(A_2)=\infty$,  then $B_1\not\in \Phi(K_1)$, and so
$\beta(A_2)=\alpha(B_1)=\infty$ since that $B_1$ is right
invertible. Otherwise, if $\beta(A_2)<\infty$, then both $A_2$ and
$B_1$ are Fredholm. Consequently,

$$0=\makebox{ind}(\left(\begin{array}{cc} A_2&C_{21}\\0&B_1\\
 \end{array}\right))=\makebox{ind}(A_2)+\makebox{ind}(B_1)=
-\beta(A_2)+\alpha(B_1),$$ that is,  $\beta(A_2)=\alpha(B_1).$
Therefore, $s.\_mul(A)=b.s.\_mul(B).$

$(2)\Rightarrow(1)$. Suppose that $A\in \Phi_{ab}(H)$, $B\in
\Phi_{sb}(K)$ and that $s.\_mul(A)=b.s.\_mul(B).$  Then from
Proposition 1.7 we have that  $A$ can be decomposed into the
following form with respect to some orthogonal decomposition
$H=H_1\oplus H_2$

$$A=\left(\begin{array}{cc}A_1&A_{12}\\0&A_2\end{array}\right),$$

\noindent where and $\dim(H_1)<\infty$, $A_1$ is nilpotent, and
$A_2$ is a left invertible operator. By Proposition 1.8, $B\in B(K)$
can be decomposed into the following form with respect to some
orthogonal decomposition $K=K_1\oplus K_2$

$$B=\left(\begin{array}{cc}B_1&B_{12}\\0&B_2\end{array}\right),$$

\noindent where $\dim(K_2)<\infty$,  $B_1$ is a right invertible
operator, and $B_2$ is nilpotent. Moreover, $s.\_mul(A)=\beta(A_2)$
and $b.s.\_mul(B)=\alpha(B_1)$. Since the assumption that
$s.\_mul(A)=b.s.\_mul(B)$, $\alpha(B_1)=\beta(A_2)$. Thus, we
conclude from Theorem 1.5 that  there exists some operator $C_{12}\in
B(K_1, H_2)$ such that $\left(\begin{array}{cc} A_2&C_{21}\\0&B_1\\
 \end{array}\right)$ is invertible. Define $C\in B(K,H)$ as follows:
$$C=\left(\begin{array}{cc}0&0\\C_{12}&0\end{array}\right).$$
By Lemma 1.6, it no hard to prove that $M_C\in \Phi_{b}(H\oplus K)$.

$(2)\Longleftrightarrow (3)$.  For this, it is sufficient to prove that if

$A\in \Phi_{ab}(H)$ and $B\in \Phi_{sb}(K)$, then
$$\alpha(A)+ \alpha(B)= \beta(A)+\beta(B)\makebox{\,\,if and only if } b.s.\_mul(B)=s\_mul(A),$$
which follows from Propositions 1.7 and 1.8 immediately. This
completes the proof.

\vskip 0.1in

In [1], Cao has proved the equivalence of $(1)$ and $(3)$ of
Theorem 2.7 by a different method, which seems to be more complicated.

\vskip 0.1in

The next corollary immediately follows from Theorem 2.7.

\vskip 0.2in \noindent {\bf Corollary 2.8.} For any given $A\in
B(H)$ and $ B\in B(K)$, we have

$\begin{array}{ll}\bigcap_{C\in G(K,\,H)}\sigma_{b}(M_{C})&
=\sigma_{ab} (A)\cup \sigma_{sb}(B)\cup\\&\,\,\,\,\{\lambda\in\Phi_{ab}
(A)\cap\Phi_{sb}(B):b.s.\_mul(B-\lambda)\not=s\_mul(A-\lambda)\}\\
&=\sigma_{ab} (A)\cup \sigma_{sb}(B)\cup\\&\,\,\,\,\{\lambda\in {\mathbb{C}}:
\alpha(A-\lambda)+ \alpha(B-\lambda)\not= \beta(A-\lambda)+
\beta(B-\lambda)\}.\\\end{array}$

\end{document}